\documentclass[oneside,12pt]{amsart} 
\usepackage{amsmath, amssymb, amsfonts, amstext, amsthm, amscd}
\usepackage[mathscr]{euscript}
\usepackage{enumerate}
\usepackage{tikz,color,soul}

\usetikzlibrary{matrix,arrows}
\usepackage[english]{babel}
\usepackage{chngcntr}
\usetikzlibrary{arrows.meta, positioning,arrows}
\newtheorem{thm}{Theorem}[section]
\newtheorem{prop}[thm]{Proposition}
\newtheorem{lemma}[thm]{Lemma}
\newtheorem{cor}[thm]{Corollary}
\theoremstyle{defn}
\newtheorem{defn}[thm]{Definition}

\newtheorem{rmk}[thm]{Remark}

\numberwithin{equation}{section}

\newcommand{\C}{\mathbb C}

\newcommand{\R}{\mathbb R}
\newcommand{\Z}{\mathbb Z}

\begin{document}

\title{Equivariant $K$-Ring of quasitoric manifolds}

\author{Jyoti Dasgupta}
\address{Department of Mathematics, Indian Institute of Technology-Madras, Chennai, India}
\email{jdasgupta.maths@gmail.com}

\author{Bivas Khan}
\address{Department of Mathematics, Indian Institute of Technology-Madras, Chennai, India}
\email{bivaskhan10@gmail.com}

\author{V. Uma}
\address{Department of Mathematics, Indian Institute of Technology-Madras, Chennai, India}
\email{vuma@iitm.ac.in}

\begin{abstract}
  Let $X(Q,\Lambda)$ be a quasitoric manifold associated to a simple
  convex polytope $Q$ and characteristic function $\Lambda$. Let
  $T\cong (\mathbb{S}^1)^n$ denote the compact $n$-torus acting on
  $X=X(Q,\Lambda)$. The main aim of this article is to give a
  presentation of the $T$-equivariant $K$-ring of $X$, as a
  Stanley-Reisner ring over $K^*(pt)$. We also derive the
  presentation for the ordinary $K$-ring of $X$.
\end{abstract}

\maketitle

\section{Introduction}
Let $X=X(Q,\Lambda)$ be a quasitoric manifold associated to a simple
convex polytope $Q$ and characteristic function
$\Lambda:\mathcal{F}\longrightarrow \mathbb{Z}^n$. Here
$\mathcal{F}:=\{Q_1,\ldots, Q_d\}$ denote the facets of $Q$. Recall
that the ordinary $K$-ring of a point is
$K^*(pt)=\mathbb{Z}[z,z^{-1}]$. Here $z$ denotes the Bott periodicity
element having cohomological dimension $-2$. Define the ring
\begin{equation}\label{kfr}
\mathcal{K}(Q,\Lambda):=\frac{K^*(pt)[y_1^{\pm1},\ldots,y_d^{\pm1}]}{J}\end{equation}
where $J$ is the ideal generated by elements of the form
\[ (1-y_{i_1})\cdots(1-y_{i_r}) ~ \mbox{whenever}
  ~Q_{i_1}\cap\cdots\cap Q_{i_r}=\emptyset.\] We call
$\mathcal{K}(Q,\Lambda)$ as the {\em $K$-theroretic face ring of the
  polytope $Q$}. We refer to \cite[Section 2]{segal} for the
definition of equivariant $K$-theory. In particular, $K^0_{T}(X)$
denotes the Grothendieck ring of $T$-equivariant complex vector
bundles on $X$. There is a canonical ring structure on $K^0_{T}(X)$
induced by the operations of tensor product and direct sum of
$T$-equivariant vector bundles on $X$. This ring structure extends to
$K_{T}^*(X)$ making it into a graded ring.  In particular,
$K_{T}^0(pt)=R(T)$ is the Grothendieck ring of finite dimensional
complex representations of $T$ and
\(\displaystyle K_{T}^*(pt)=K_{T}^0(pt)\bigotimes_{\mathbb{Z}}
K^*(pt)=R(T)[z,z^{-1}]\). The ring $K_{T}^*(X)$ gets the structure of
$K_{T}^*(pt)$-algebra via the equivariant canonical projection
$X\longrightarrow pt$.

In our main theorem Theorem \ref{main} we prove that we have a
canonical isomorphism $\mathcal{K}(Q,\Lambda) \simeq K_{T}^*(X)$ which
takes $y_{i}$ to $[L_{i}]$ where $L_{i}$ is a canonical
$T$-equivariant complex line bundle on $X$ corresponding to the facet
$Q_{i}$ (see \cite[Lemma $2.3$]{torusmanifold}, \cite[Section
6.2]{davisjanu} and \cite[Section 3]{SU}). We note that Theorem
\ref{main} is analogous to the corresponding result for algebraic
equivariant $K$-ring of smooth toric varieties in \cite[Theorem
6.4]{vezzosi2003higher} and \cite[Proposition
2.2]{baggio2007equivariant}. This is also a generalization of the
description of topological equivariant $K$-ring of toric varieties
associated to a smooth polytopal fan in \cite[Corollary 7.9]{HHRW}
(See \cite[p. 419 $(C')$]{davisjanu}).

Furthermore, since $X$ is weakly equivariantly formal (see
\cite[Definition 4.1]{HaLa}), we recover the presentation of $K^*(X)$
in Corollary \ref{maincor}. The relative version of this result was
proved using different techniques in \cite[Theorem 1.2]{SU}.

\section{Notations and Preliminaries}
We recall some notations and preliminary facts from \cite{davisjanu}
and \cite{SU}.  An $n$-dimensional convex polytope $Q$ is said to be
simple if precisely $n$ facets (codimension $1$ faces) meet at each
vertex. Let $T \cong \left( S^1 \right)^n$ denote the $n$-dimensional
compact torus and let $Q \subset \R^n$ be a simple convex polytope of
dimension $n$. A (smooth) $T$-manifold $X$ with a locally standard
action of $T$ and projection $\pi: X \rightarrow X/T \cong Q$, is
called a \emph{$T$-quasitoric manifold over $Q$}. Here `local
standardness' means that $X$ has a covering by $T$-invariant open sets
$U$ such that $U$ is weakly equivariantly diffeomorphic to an open
subset $U'\subset \C^n$ invariant under the standard $T$-action on
$\C^n$.  The latter means that there is an automorphism
$\theta :T \rightarrow T$ and a diffeomorphism $f:U \rightarrow U'$
such that $f(ty)=\theta(t) f(y)$ for all $t \in T$, $y \in U$. Let $F$
be a face of $Q$. Then for any $x \in \pi^{-1} (\text{int } F)$, the
isotropy group at $x$ is independent of the choice of $x$, we denote
it by $T_F$. Also note that the rank of $T_F$ is same as the
codimension of $F$.

Let $\mathcal{F}:=\{ Q_1, \ldots, Q_d \} $ denote the set of facets of
$Q$. For each facet $Q_j$, the rank-$1$ subgroup $T_{Q_j}$ is
determined by a primitive vector
${\lambda}_j \in \Z^n\cong \text{ Hom}(S^1, T)$ which is unique upto
sign. Consider $V_j = \pi^{-1}(Q_j)$, then each $V_j$ is
orientable. The sign of ${\lambda}_j$ is determined by choosing an
omniorientation on $X$, i.e. orientation on $X$ as well as on each
$V_j$ for $1 \leq j \leq d$. Then we can define the {\it
  characteristic map}
$\Lambda : \mathcal{F} \rightarrow \Z^n \cong \text{Hom }(S^1, T)$
which sends $Q_j$ to ${\lambda}_j$. The map $\Lambda$ also satisfies
the following condition:
\begin{multline}\label{smooth condition}
\text{If }Q_{i_1}\cap \cdots \cap Q_{i_k}\text{ is a face of }P,\text{ then }
\Lambda(Q_{i_1}), \ldots , \Lambda (Q_{i_k}) \\
\text{ is a part of a basis for the integral lattice Hom}(S^1, T)\cong \Z^n.
\end{multline}
 
Conversely, given any $n$-dimensional simple convex polytope $Q$ and a characteristic map $\Lambda$, there exists an omnioriented quasitoric manifold $X$ over $Q$ with characteristic map $\Lambda$. The omnioriented quasitoric manifold $X$ is determined upto equivalence over $Q$ by its characteristic function.\\

Let $X$ be an omnioriented quasitoric manifold over $Q$. Realizing $Q$
as a convex polytope in $\R^n$, we can choose a vector ${w} \in \R^n$
such that ${w}$ is tangent to no proper face of $Q$. Then we have the
linear function $h :\R^n \rightarrow \R$ given by
$h(x)= \langle {x}, {w} \rangle$ where $\langle , \rangle$ is the
Euclidean inner product on $\R^n$. We can view $h$ as a height
function on $Q$, i.e, $h$ is injective when restricted to any edge of
$Q$. The map $h$ induces an ordering on the set of vertices
$\mathcal{V}$ of $Q$, where $v<v'$ in $\mathcal{V}$ if $h(v)
<h(v')$. Also this induces an orientation on the edges which makes the
$1$-skeleton of $Q$ into a directed graph by directing the edge
$E_{v,v'}$ joining $v$ and $v'$ by $v\rightarrow v'$ whenever
$v<v'$. Since $Q$ is simple, exactly $n$ edges meet at each
vertex. For each vertex $v$ of $Q$, let $\text{ind}(v)$ be the number
of inward pointing edges at $v$ and let $Q_v$ be the face of $Q$
spanned by the incoming edges at $v$. We denote by $\widehat{Q}_v$ to
be the subset obtained from $Q_v$ by removing all faces not containing
$v$ and let $Z_v= \pi^{-1}(\widehat{Q}_v)$. Then $\widehat{Q}_v$ is
homeomorphic to $\R_{\geq}^{ind(v)}$ and $Z_v$ can be identified with
$\C^{ind(v)}$, i.e. it is a cell of dimension $2 \cdot ind(v)$. Also
note that the closure of the cell $Z_v$ is the quasitoric submanifold
$\pi^{-1}(Q_v)$. The union of subsets $Z_v$ over $v\in \mathcal{V}$
defines a cell decomposition of $X$ with cells only in even
dimensions.

We fix some further notations below.

Let $m=|\mathcal{V}|$ and $v_1< \cdots < v_m$ be the ordering of
$\mathcal{V}$ associated to $h(v_i)<h(v_{i+1})$ for $1\leq i\leq m-1$.

Let $Z_i:=\pi^{-1}(\widehat{Q}_{v_i})$ for $1\leq i\leq m$. Also let
$x_i:=\pi^{-1}(v_i)\in Z_i$ denote the $T$-fixed points in $X$ and
$k_i:=\text{ind}(v_i)$ for $1\leq i\leq m$. In particular, $k_1=0$ and
$k_m=n$.

Let \(\displaystyle X_i:=\bigcup_{j\leq i} Z_j\). Thus
$X_1=\{x_1\}\subseteq X_2\subseteq \cdots\subseteq X_m=\{X\}$ defines
a stratification of $X$ by $T$-invariant submanifolds such that
$X_i\setminus X_{i-1}=Z_i\cong \mathbb{C}^{k_i}$.

\section{GKM theory on quasitoric manifolds}
We begin this section by recalling the GKM theory from \cite[Section
3]{HHH}. We shall then verify that these results can be applied to a
quasitoric manifold and hence give a precise description of its
topological equivariant $K$-ring.

Let $Y$ be a $G$-space for a compact Lie group $G$ equipped with a $G$-invariant stratification
\[Y_1\subseteq Y_2\subseteq\cdots\subseteq Y_m=Y.\] Every
$Y_i\setminus Y_{i-1}$ has a $G$-invariant subspace $F_i$ whose
neighbourhood carries the structure of the total space of a
$G$-equivariant vector bundle $\rho_i=({V}_i,\pi_i,F_i)$.

Let $K_{G}^*(Y)$ denote the $G$-equivariant $K$-ring of the $G$-space
$Y$.

\begin{enumerate}

\item[A1] Each subquotient $Y_i/Y_{i-1}$ is homeomorphic to the Thom
  space $Th(\rho_i)$ with attaching maps
  $\phi_i:S(\rho_i)\longrightarrow Y_{i-1}$.

\item[A2] Every $\rho_i$ admits a $G$-equivariant direct sum
  decomposition \(\displaystyle\bigoplus_{j<i} \rho_{ij}\) into $G$-equivariant
  subbundles $\rho_{ij}=(V_{ij},\pi_{ij}, F_i)$. We allow the case $V_{ij}=0$.

\item[A3] There exist $G$-equivariant maps  $f_{ij}:F_i
  \longrightarrow F_j$ such that the restrictions $f_{ij}\circ
  \pi_{ij}\mid_{S(\rho_{ij})}$ and 
$\phi_i\mid_{S(\rho_{ij})}$ agree for every $j<i$. 

\item[A4]The equivariant Euler classes $e_G(\rho_{ij})$ for $j<i$, are not
   zero divisors and are pairwise relatively prime in $K^*_G(F_i)$. 

\end{enumerate}

We now recall \cite[Theorem 3.1]{HHH} which gives a precise description of the
generalized $G$-equivariant $K$-theory of $Y$.

\begin{thm}\label{gkm}
Let $Y$ be a $G$-space satisfying the four assumptions A1 to A4. Then
the restriction map
\[\iota^*: K^*_{G}(Y)\longrightarrow \prod_{i=1}^m K_{G}^*(F_i)\] is
monic and its image $\Gamma_{Y}$ can be described as
\[\{(a_i)\in \prod_{i=1}^m K_G^*(F_i):
e_{G}(\rho_{ij})~\mid~a_i-f_{ij}^*(a_j)~\forall~ j<i)\}.\]
\end{thm}

We first show below that the GKM theory of \cite{HHH} described above
can be applied to the $T$-space $X$. Further, we shall use this to
give an explicit description of $K_T^{*}(X)$.

\begin{rmk}\label{generalized}
  Note that Theorem \ref{gkm} has been proved in \cite{HHH} for the
  generalized $G$-equivariant cohomology theory $E_{G}^*(Y)$.
\end{rmk}

\begin{prop}\label{divass}
  A quasitoric manifold $X(Q,\Lambda)$ with the above
  $T$-invariant stratification satisfies assumptions A1 to
  A4 listed above. 
\end{prop}

{\bf Proof:} Since $Z_i\simeq \mathbb{C}^{k_i}$ are $T$-stable we have
the associated $T$-representation
\[\rho_i=(V_i:=\mathbb{C}^{k_i},\pi_i,x_i),\] which can
alternately be viewed as a $T$-equivariant complex vector bundle over
the $T$-fixed point $x_i$ for $1\leq i\leq m$.

Let $v_{i_1}, v_{i_2},\ldots, v_{i_{k_i}}\in \mathcal{V}$ be such that
there are directed edges $E_{i_1},\ldots, E_{i_{k_i}}$ respectively
from $v_{i_1}, v_{i_2},\ldots, v_{i_{k_i}}$ to $v_i$ in the directed
graph on the $1$-skeleton of $Q$ associated to the height function
$h$. Since $Q$ is simple, on the curve $C_{i_j}:=\pi^{-1}(E_{i_j})$,
the torus $T$ acts via a character $\chi_{i_j}:T\longrightarrow
S^1$. Let $E_{i_j}=Q_{j_1}\cap\cdots\cap Q_{j_{n-1}}$. The
$n\times (n-1)$-matrix
$[\Lambda(Q_{j_1}),\cdots ,\Lambda(Q_{j_{n-1}})]$, defines a
$\mathbb{Z}$-linear map
$\psi_j:\mathbb{Z}^n\longrightarrow \mathbb{Z}^{n-1}$ of
$\mathbb{Q}$-rank $(n-1)$, via multiplication from the right by
viewing elements of $\mathbb{Z}^n$ and $\mathbb{Z}^{n-1}$ as
$1\times n$ and $1\times (n-1)$ integer matrices respectively.  The
character $\chi_{i_j}$ is determined by a primitive vector
$u_{i_j}\in \mathbb{Z}^n$ generating the kernel of $\psi_j$ whose sign
is determined by the omniorientation.

Let $N(i):=\{i_1,\ldots, i_{k_i}\}$. Since
$\mathbb{R}^{k_i}_{\geq 0}\simeq \widehat{Q}_{v_i}\subset Q$ is
generated by $E_{i_j}\setminus v_j$ for $1\leq j\leq k_i$ it follows
that
\(\displaystyle\rho_i=\bigoplus_{j=1}^{k_i} \mathbb{C}_{\chi_{i_j}}\)
where $\mathbb{C}_{\chi_{i_j}}$ is the one dimensional complex
representation $\rho_{i_j}$ determined by the character
$\chi_{i_j}$. Thus $\rho_i$ is a direct sum of $1$-dimensional
$T$-representations $\rho_{ij}$ where $T$ acts by the character
$\chi_{i_j}$. Thus if we let $V_{il}:=\mathbb{C}_{\chi_{i_j}}$ for
$l=i_j\in N(i)$ and $V_{il}=0$ for $l<i$ and $l\notin N(i)$ then this
verifies assumption (2).

Consider the neighbourhood $W'_i=\widehat{Q}_{v_i}\cap D$ of $v_i$ in
$Q$ where $D$ is a closed disc in $\mathbb{R}^n$ with centre $v_i$ not
containing any other vertex of $Q$. Now $Q_{v_i}$ is the smallest face
of $Q$ containing the edges $E_{i_j}$ for $1\leq j\leq k_i$.  Note that the
$Link(v_i)$ in $Q_{v_i}$ is
\(\displaystyle\bigcup_{l<i} \widehat{Q}_{v_l}\) and
$Q_{v_i}=Star(v_i)=Link(v_i)\star v_i$. Thus it follows from
polyhedral geometry that for every
\(\displaystyle p\in \bigcup_{l<i} \widehat{Q}_{v_l}\), the line
segment in $Q$ joining $p$ and $v_i$ meets $\partial (W'_i)$ at a
unique point $y_i$. Moreover, $y_i$ determines $p=p(y_i)$ uniquely and
vice versa. This gives a bijective correspondence between
$\partial(W'_i)$ and
\(\displaystyle\bigcup_{l<i} \widehat{Q}_{v_l}\). Note that if
$y_i$ belongs to the relative interior of a face $F$ of $Q_{v_i}$ then
$p(y_i)\in F$. Hence $t\cdot y_i=y_i$ implies $t\cdot
p(y_i)=p(y_i)$. Thus we get an attaching map
\(\displaystyle g_i:\pi^{-1}(\partial W'_i)\longrightarrow
\bigcup_{l<i} Z_l\subseteq X_{i-1}\) that maps $[t,y_i]$ to
$[t,p(y_i)]$. Moreover, under the identification of $Z_i$ with the
complex representation $\rho_i$,
$W_i:=\pi^{-1}(W_i')\subseteq Z_i\subseteq X_i$ is identified with the
disc bundle $D(\rho_{i})$ and
$\partial(W_i):=\pi^{-1}(\partial(W_i'))$ with the sphere bundle
$S(\rho_i)$ of the representation space. This induces the following
homeomorphisms
$X_i/X_{i-1}\simeq W_i/\partial(W_i)\simeq
D(\rho_i)/S(\rho_i)=Th(\rho_i)$, where $Z_i\subseteq X_i$ is
identified with the interior of $W_i$. This verifies assumption (1).

Let $f_{ij}:x_i\longrightarrow x_{i_j}$
denote the constant map for every $1\leq j\leq k_i$. Further, if
$\widehat{E}_{i_j}:= (E_{i_j}\setminus v_{i_j})\subseteq Q_{v_i}$, then
$\pi^{-1}(\widehat{E}_{i_j})$ can be identified with the one dimensional
sub-representation $\rho_{i_j}$ of $\rho_i$. Let $S(\rho_{i_j})$ be
the circle bundle associated with $\rho_{i_j}$. Let $w_{i_j}$ denote
the point where $E_{i_j}$ meets $\partial(W'_i)$. Then the attaching
map $g_i$ sends $\pi^{-1}(w_{i_j})$ in $S(\rho_{i_j})$ to
$x_{i_j}$. Further, $\pi^{-1}(w_{i_j})\in W_i$ is mapped to $x_i$
under the canonical projection of the representation $\rho_{i_j}$ and
further mapped to $x_{i_j}$ after composition with $f_{ij}$. It
follows that the restrictions to $S(\rho_{i_j})$ of $g_i$ and
the composition of the projection of $\rho_{i_j}$ with $f_{ij}$ agree
for every $1\leq j\leq k_i$. This verifies assumption (3).

Note that for $1\leq j\neq j'\leq k_i$, $E_{i_j}$ and $E_{i_{j'}}$
are distinct edges incident at $v_i$ in the simple polytope $Q$. Thus
if $v_i=Q_{i_1}\cap\cdots\cap Q_{i_n}$ then clearly
$E_{i_{j'}}\subseteq Q_{i_r}$ where $i_r\in \{i_1,\ldots,i_n\}\setminus \{j_1,\ldots,j_{n-1}\}$.  Since
$\Lambda(Q_{i_1}),\cdots ,\Lambda(Q_{i_n})$ form a $\mathbb{Q}$ basis, and $u_{i_j}$ is orthogonal to $\Lambda(Q_{j_1}),\ldots \Lambda(Q_{j_{n-1}})$ and  $u_{i_{j'}}$ is orthogonal to $\Lambda(Q_{i_r})$ it
follows that $u_{i_j}$ and $u_{i_{j'}}$ are linearly independent.

Now, $u_{i_j}$ is primitive and non-zero for every $1\leq j\leq k_i$,
the $K$-theoretic equivariant Euler class
$e^{T}(V_{il})=1-e^{-u_{i_j}}$ for $l=i_j\in N(i)$, is a non-zero
divisor in the ring $K^0_T(x_i)=RT$, which is an integral domain. Here
$RT$ denotes the ring of isomorphism classes of finite dimensional
complex representations of $T$. Recall that
$RT=\mathbb{Z}[e^u:u\in M]$ where $M=Hom(T,\mathbb{C}^*)$ is the
character group of $T$.  Moreover, since $u_{i_j}$ and $u_{i_{j'}}$
are linearly independent, it follows that $1-e^{-u_{i_j}}$ and
$1-e^{-u_{i_{j'}}}$ are relatively prime in the U.F.D $RT$ (see
\cite[Lemma 5.2]{HHH}).  Also if $l<i$ and $l\notin N(i)$, then
$V_{il}=0$ so that $e^{T}(V_{il})=1$. This verifies assumption
(4). $\Box$

\begin{rmk}\label{composeconstant}
  Note that we can define the constant map
  $f_{ij}:x_i\longrightarrow x_j$ between any two $T$-fixed points of
  $X=X(Q,\lambda)$, which satisfy $f_{ik}=f_{jk}\circ f_{ij}$ for
  $1\leq i,j,k \leq m$. Thus we have the pull-back map of equivariant
  $K$-theory $f_{ij}^*: K^*_{T}(x_j)\longrightarrow K^*_{T}(x_i)$
  which satisfies $f_{ik}^*=f_{ij}^*\circ f_{jk}^*$ for
  $1\leq i,j,k\leq m$. The canonical maps $s_i:pt\longrightarrow x_i$
  induce isomorphisms $s_i^*:K_{T}^*(x_i)\longrightarrow K_{T}^*(pt)$
  for $1\leq i\leq m$. Under identifications of $K_{T}^*(x_i)$ with
  $K_{T}^*(pt)$ via $s_{i}^*$ $1\leq i\leq m$,
  $f_{ij}^*=(s_{i}^*)^{-1}\circ s_j^*$ correspond to
  $id\mid_{K_{T}^*(pt)}$ for $1\leq i,j\leq m$. \end{rmk}

\section{Equivariant K-theory of a quasitoric manifold}

\begin{prop}\label{eqkt}
Let $X=X(Q,\Lambda)$ be a  quasitoric manifold. The $T$-equivariant ring
$K_T^*(X)$ of $X$ is isomorphic to the subring 
\[\Gamma_{X}=\{(a_i)~:~1-e^{-u_{i_j}}~\mid~
  a_i-a_{i_j}~\forall~i_j\in N(i)\}\subseteq {K_{T}^*(pt)}^m.\]
\end{prop} {\bf Proof:} By Proposition \ref{divass} and Theorem
\ref{gkm} above it follows that $K^*_{T}(X)$ is isomorphic to the
subring
\begin{equation}\label{localization1}
  \{(a_i)~:~e^T(V_{il})~\mid~a_i-f_{il}^*(a_l)~\mbox{for every}~ l<i\}\subseteq \prod_{i=1}^m K_{T}^*(x_i).\end{equation} 
Now, from the proof
of Proposition \ref{divass} we have $V_{il}=0$ for $l<i$ and
$l\notin N(i)$ and $V_{il}=\rho_{i_j}$ for $l=i_j\in N(i)$. Thus
$e^T(V_{il})=1$ for $l<i$ and $l\notin N(i)$ and
$e^{T}(V_{il})=1-e^{-u_{i_j}}$ for $l=i_j\in N(i)$. Hence we can rewrite
(\ref{localization1}) as
\begin{equation}\label{localization2}\{(a_i)~:~1-e^{-u_{i_j}}~\mid~
  a_i-f^*_{ii_j}(a_{i_j})~\mbox{for every}~i_j\in N(i) \subseteq \prod_{i=1}^m K_{T}^*(x_i).\end{equation} By Remark
\ref{composeconstant}, (\ref{localization2}) can be further identified with
$\Gamma_{X}$.  Hence the proposition.  $\Box$

\begin{rmk}\label{edges}
  Since for every edge $E_{il}$ either $i\in N(l)$ or $l\in N(i)$
  depending on whether $h(v_i)<h(v_l)$ or vice versa. Thus the
  subring $\Gamma_{X}$ can equivalently be defined as
  \begin{equation}\label{gkm1}\Gamma_{X}=\{(a_i)~:~e_{T}({V_{il}})~\mid~a_i-a_l~\mbox{for
      every edge}~ E_{il} ~\mbox{in}~Q\}\subseteq K_{T}^*(pt)^m.\end{equation} Here $V_{il}$ is the
  $T$-representation corresponding to the primitive character
  $u_{il}\in \mathbb{Z}^n$ orthogonal to $\Lambda(Q_{(il)_k})$, $1\leq k\leq
  n-1$ where \(\displaystyle E_{il}=\bigcap_{k=1}^{n-1}
  Q_{(il)_k}\). Thus $u_{il}=u_{i_j}$ or $u_{l_t}$ depending on
  whether $l=i_j\in N(i)$ or $i=l_t\in N(l)$.
\end{rmk}

\begin{rmk} Proposition \ref{divass} and Proposition \ref{eqkt} have
  been proved for the more general setting of a divisive toric
  orbifold in \cite[Proposition 4.2 and Proposition
  5.1]{saruma}. \end{rmk}

Recall that $RT=\Z[M]=\Z[e^u ~:~ u \in M]$ where $M=X^*(T)$
is the character group of $T$. Let
\[R:=K_{T}^*(pt)=K^*(pt)\otimes_{\mathbb{Z}}RT= K^*(pt)[e^u : u\in
  M].\] For any face $F=Q_{i_1}\cap \cdots \cap Q_{i_l}$ of $Q$, let
\[R_{F}:=R/{I_{F}}\] where $I_{F}$ is the ideal in $R$ generated by
the elements $1-e^u$ for
$u\in \langle\Lambda(Q_{i_1} ),\ldots,
\Lambda(Q_{i_l})\rangle^{\perp}$.  Clearly, for $F' \prec F$, we have
$I_{F'} \subseteq I_F$. This induces a natural projection map
\begin{equation}\label{np} f_{F',F}:R_{F'} \rightarrow R_F\end{equation} which
  satisfies \begin{equation}\label{composeres} f_{F',F} \circ
    f_{F'',F'}=f_{F'',F}\end{equation} whenever $F'' \prec F' \prec F$.
 
  In particular, for any $v \in \mathcal{V}$ we see that \(R_v=R.\)
  Now, for any $w \in \mathcal{V}$, denote by $v \vee w$ the minimal
  face of $Q$ containing both $v$ and $w$. If $v \vee w=Q$, then
  clearly $R_{v \vee w}=\Z$ and the natural projection map
  $f_{v,v\vee w}:R_v \rightarrow R_{v\vee w}$ is the augmentation
  map. On the other hand when $v \vee w$ is a proper face, write
  $v \vee w=Q_{i_1} \cap \cdots \cap Q_{i_k}$. Then
  $R_{v \vee w}:=\frac{R}{I_{v \vee w}}$ where $I_{v \vee w}$ is the
  ideal in $R$ generated by the elements $1-e^u$ for
  $u\in \langle\Lambda(Q_{i_1} ),\ldots,
  \Lambda(Q_{i_k})\rangle^{\perp}$. Then we have the canonical
  projection map $R_v \rightarrow R_{v \vee w}$ which sends $e^{u}$ to
  $1$ for
  $u\in \langle
  \Lambda(Q_{i_1}),\ldots,\Lambda(Q_{i_k})\rangle^{\perp}$. If
  $E_{il}= v_i \vee v_l$ is an edge, then $I_{E_{il}}$ is the
  principal ideal generated by $1-e^{u_{il}}$ where $u_{il}$ is the
  primitive character in $\mathbb{Z}^n$ orthogonal to
  $\Lambda(Q_{(il)_1}),\ldots, \Lambda(Q_{(il)_{n-1}})$, where
  \(\displaystyle E_{il}= Q_{(il)_1}\cap\cdots \cap
  Q_{(il)_{n-1}}\) (see Remark \ref{edges}).

 \begin{lemma} The ring $\Gamma_X$  in Proposition \ref{eqkt} can be rewritten as
   \[\{(a_i) \in R^m : f_{v_i,E_{il}}(a_i)=f_{v_l,E_{il}}(a_l) ~
     ~\text{whenever there is an edge} ~E_{il} ~~\text{in}~ Q
     ~\text{joining} ~v_i~\text{and} ~v_l \}.\]\end{lemma} {\bf
   Proof:} Since $I_{E_{il}}$ is generated by $1-e^{u_{il}}$, by
 Proposition \ref{eqkt} and Remark \ref{edges} we have
 \[\Gamma_X=\{(a_i) \in R^m
   ~\mid~ a_i-a_l \in I_{E_{il}} ~\text{whenever there is an edge}~E_{il}
   ~\text{joining} ~v_i ~\text{and}~v_l\}.\] By (\ref{np}) this further implies that
\[\Gamma_{X}=\{ (a_i) \in R^m ~\mid~
  f_{v_i,E_{il}}(a_i)=f_{v_l,E_{il}}(a_l) ~ \text{whenever there is an edge}
  ~E_{il} ~\text{joining} ~v_i ~\text{and}~v_l\}.\] Hence the lemma.
$\Box$

Let
\(\mathcal{W}_{X}:= \{ (a_i) \in R^m ~:~ f_{v_i,F}(a_i)=
f_{v_l,F}(a_l) ~\text{ for } i \neq l ~\text{and}~ F=v_i \vee v_l\}.\)
Clearly by definition we have an inclusion
$\mathcal{W}_{X} \subseteq \Gamma_X$ as subrings of $R^m$. The
following lemma shows that the inclusion is in fact an equality.

\begin{lemma}
We have an equality $\Gamma_X=\mathcal{W}_{X}$ as subrings of $R^m$.
\end{lemma} 
{\bf Proof:} Let $(a_i) \in \Gamma_X$. For
$v_i, v_l \in \mathcal{V}$, let
\[v_i=v_{i^0} < v_{i^1} < \cdots < v_{i^k}=v_l\] be the minimal
sequence of vertices in $\mathcal{V}$ such that any there is an edge
joining any two consecutive vertices. Then
$E_{i^p, i^{p+1}} \prec F=v_i \vee v_l$ for $p=0, \ldots, k-1$. Now,
we need to show that $ f_{v_i,F}(a_i)= f_{v_l,F}(a_l)$. Since
$(a_i)\in \Gamma_X$, for every $0\leq p\leq {k-1}$ we have,
\[f_{v_{i^{p}},E_{i^p,
      i^{p+1}}}(a_{i^p})=
  f_{v_{i^{p+1}},E_{i^p,
      i^{p+1}}}(a_{i^{p+1}}).\] Hence 
\[ f_{E_{i^p, i^{p+1}} ,F}\circ f_{v_{i^p},E_{i^p, i^{p+1}}}(a_{i^p})=
  f_{E_{i^p, i^{p+1}} , F} \circ f_{v_{i^{p+1}},E_{i^p,
      i^{p+1}}}(a_{i^{p+1}}) .\] This implies by (\ref{composeres})
that $f_{v_{i^p},F}(a_{i^p})= f_{v_{i^{p+1}},F}(a_{i^{p+1}})$ for
$p=0, \ldots k-1$.  Thus we get
$f_{v_i,F}(a_i)= f_{v_l,F}(a_l)$.$\hfill \Box$

\begin{defn}\label{14} The $K$-theoretic face ring of the polytope $Q$
  is defined to be
 \(\displaystyle \mathcal{K}(Q, \Lambda):=\frac {K^*(pt)[y_1^{\pm 1}, \ldots,
    y_d^{\pm 1}]}{J} \) where $J$ is the ideal generated by
  elements of the form
\begin{equation}
\left(1- y_{i_1} \right)\cdots\left(1- y_{i_r} \right)~ \text{whenever} ~Q_{i_1} \cap \cdots \cap Q_{i_r} = \emptyset.
\end{equation}
 \end{defn}
 Note that one has a monomorphism of rings
 $R \stackrel{\theta}{\rightarrow}\mathcal{K}(Q, \Lambda)$ defined by
 \(\displaystyle e^u \mapsto \prod_{1 \leq i \leq d} y_{i}^{\langle u,
   \lambda_i \rangle}\), $u \in M$, which gives an $R$-algebra
 structure on $\mathcal{K}(Q, \Lambda)$. Also \(\displaystyle R^m\)
 has a canonical $R$-algebra structure via the diagonal embedding
 $\Delta$.

 For $1\leq i\leq d$ such that $v_j\in Q_i$, using (\ref{smooth
   condition}), we define $\mu_i^j\in M$ by
 $\langle\mu_i^j,\lambda_l \rangle=\delta_{i,l} ~(\mbox{Kronecker
   delta})$ for every $1\leq l\leq d$ such that $v_j\in Q_l$. For
 $1\leq i\leq d$ such that $v_j\notin Q_i$, we define $\mu^j_i=0$.

 Define the map
 \(\displaystyle \phi:K^*(pt) [y_{1}^{\pm 1}, \ldots, y_{d}^{\pm 1}]
 \rightarrow R^m\) given by $y_i \mapsto r_i$ where $r_i$ is given by
\begin{equation}\label{r_i} (r_{i})_{j} := 
  e^{{\mu^j_i}} \end{equation}

\begin{thm}\label{7}  The map $\phi$ induces a injective $R$-algebra
  homomorphism
  \[\bar{\phi}: \mathcal{K}(Q,\Lambda) \hookrightarrow  R^m\] with
  image the $R$-subalgebra $\mathcal{W}_{X}$.
\end{thm}
\noindent
{\bf Proof:} First let us note that  $r_i \in \mathcal{W}_{X}$ for $1 \leq i \leq d$.
Let $v_p, v_q \in \mathcal{V}$ be distinct. If $v_p \vee v_q=Q$ nothing to prove. Otherwise write $F=v_p \vee v_q =Q_{i_1} \cap \cdots \cap Q_{i_l}$, where $v_p=Q_{i_1} \cap \cdots \cap Q_{i_n}$ and $v_q=Q_{i_1} \cap \cdots \cap Q_{i_l} \cap Q_{j_{l+1}} \cap \cdots \cap Q_{j_n}$. Now consider the following cases:
\begin{enumerate}
\item $v_p, v_q \notin Q_i$: Then $(r_{i})_p=1=(r_{i})_q$, hence $f_{v_p,F}((r_{i})_p) =f_{v_q,F}((r_{i})_q)$.
\item $v_p \notin Q_i$ and $v_q \in Q_i$: Then $(r_{i})_p=1$ and
  $(r_{i})_q=e^{\mu_i^q}$. Note that, $f_{v_q,F}$ maps $e^{\mu_i^q}$ to $1$, since
  $\mu_i^q \in \langle \Lambda(Q_{i_1} ),\ldots, \Lambda(Q_{i_l})
  \rangle^{\perp}$. Hence we are done in this case.

\item $v_p, v_q \in Q_i$: Then $f_{v_p,F}((r_{i})_p) =e^{\mu_i^p}$ and
  $f_{v_q,F}((r_{i})_q)=e^{\mu_i^q}$. By definition
  $\mu_i^p-\mu_i^q\in \langle \Lambda(Q_{i_1} ),\ldots,
  \Lambda(Q_{i_l}) \rangle^{\perp}$. Hence it follows that
  $f_{v_p,F}((r_{i})_p) =e^{\mu_i^p}=e^{\mu_i^q}=f_{v_q,F}((r_{i})_q)$
  in $R/I_{F}$.
\end{enumerate} 
This proves that $r_i \in \mathcal{W}_{X}$ for $1 \leq i \leq d$. We show that elements of $\mathcal{W}_{X}$ can be written as Laurent polynomials in $r_i$'s with coefficients in $K^*(pt)$.

Let $\alpha=\left( \alpha_i \right) \in \mathcal{W}_{X}$. Let $v_1=Q_{i_1} \cap \cdots \cap Q_{i_n}$, then $\alpha_1 \in R$ and hence we can find a Laurent polynomial $p_1(y_{i_1}, \ldots, y_{i_n})$ with coefficients in $K^*(pt)$ such that $p_1(r_{i_1}, \ldots,  r_{i_n})_1=\alpha_1$. Let $\alpha^1:= \alpha-p_1(r_{i_1}, \ldots, r_{i_n})$. Then we see that $\alpha^1_1=0$.

Now let
$v_2=Q_{i_1} \cap \cdots \cap Q_{i_l} \cap Q_{j_{l+1}} \cap \cdots
\cap Q_{j_n}$ such that
$F=v_1 \vee v_2= Q_{i_1} \cap \cdots \cap Q_{i_l}$. Similarly as above
there is a Laurent polynomial
$p_2(y_{i_1}, \ldots, y_{i_l}, y_{j_{l+1}}, \ldots, y_{j_n})$ with
coefficients in $K^*(pt)$ such that
\begin{equation}\label{eq0}p_2(r_{i_1}, \ldots, r_{i_l}, r_{j_{l+1}}, \ldots,
r_{j_n})_2=\alpha^1_2.\end{equation} By (\ref{r_i}) we note that
\[p_2( r_{i_1}, \ldots, r_{i_l}, r_{j_{l+1}}, \ldots, r_{j_n})_1=p_2
  (e^{\mu_{i_1}^1}, \ldots, e^{\mu_{i_l}^1}, 1, \ldots, 1)\] whose
projection to $R_F$ remains unchanged,
i.e. \begin{equation}\label{eq1} p_2(r_{i_1}, \ldots, r_{i_l},
  r_{j_{l+1}},\ldots, r_{j_n})_1= f_{v_1,F}(p_2(r_{i_1}, \ldots,
  r_{i_l}, r_{j_{l+1}}, \ldots, r_{j_n})_1).\end{equation}

Since $\alpha^1\in \mathcal{W}_{X}$,
\begin{equation}\label{eq2} f_{v_2,F}(\alpha^1_2)=f_{v_1,F}(\alpha^1_1)=0.\end{equation} Further,
$p_2(r_{i_1}, \ldots, r_{i_l}, r_{j_{l+1}}, \ldots, r_{j_n}) \in
\mathcal{W}_{X}$ together with (\ref{eq2}) implies that
\begin{equation}\label{eq3}\begin{split}f_{v_1,F}( p_2( r_{i_1}, \ldots, r_{i_l}, r_{j_{l+1}}, \ldots,
    r_{j_n})_1 )&=f_{v_2,F}( p_2(r_{i_1}, \ldots, r_{i_l},
    r_{j_{l+1}}, \ldots, r_{j_n})_2)\\
    &=f_{v_2,F}(\alpha^1_2)=0.\end{split}\end{equation} Now, (\ref{eq1}) and
(\ref{eq3}) together imply
\[p_2(r_{i_1}, \ldots, r_{i_l}, r_{j_{l+1}}, \ldots, r_{j_n})_1=0.\]
Letting
$\alpha^2 := \alpha^1-p_2( r_{i_1}, \ldots, r_{i_l}, r_{j_{l+1}},
\ldots, r_{j_n})$, we have by (\ref{eq0}) that
$\alpha^2_1=0=\alpha^2_2$.  Repeating this process for
$v_3, \ldots, v_m$, where $v_t=Q_{t_1}\cap\cdots\cap Q_{t_n}$ for
$1\leq t\leq m$, we get that
$\alpha^m_1=\alpha^m_2=\cdots=\alpha^m_m=0$, for
\(\displaystyle\alpha^m:=\alpha-\sum_{t=1}^m p_t(r_{t_1},\ldots,
r_{t_n})\). Thus $\alpha^m=0$, so that $\alpha$ is in the image of
$\phi$. Since $\alpha\in \mathcal{W}_{X}$ was arbitrary, $\phi$ is
surjective. It remains to show that
$\text{ker}(\phi)=\langle\left(1- y_{i_1} \right) \cdots \left(1-
  y_{i_r} \right) ~:~ Q_{i_1} \cap \cdots \cap Q_{i_r} = \emptyset
\rangle$.

For $v \in \mathcal{V}$, consider the map
$\phi_v:K^*(pt)[y_1^{\pm1},\ldots, y_d^{\pm 1}] \longrightarrow R$
which sends $y_i \mapsto (r_{i})_v$ for $1 \leq i \leq d$. We see that
$\text{ker}(\phi_v)=J_v:=\langle y_j-1: v \notin Q_j \rangle $ and
clearly
\(\displaystyle\text{ker}(\phi)= \bigcap_{v \in \mathcal{V}} J_v\). It
follows from \cite[Lemma $6.5$]{vezzosi2003higher} that
\(\displaystyle\bigcap_{v \in \mathcal{V}} J_v= \langle\left(1-
  y_{i_1} \right) \cdots \left(1- y_{i_r} \right) ~:~ Q_{i_1} \cap
\cdots \cap Q_{i_r} = \emptyset \rangle\). Hence we get the induced
ring homomorphism $\bar{\phi}$ from
$\mathcal{K}(Q,\Lambda)\longrightarrow \mathcal{W}_{X}\subseteq
R^m$. It follows by (\ref{smooth condition}) that for each
$1\leq j\leq m$, $\{\mu^j_i~:~1\leq i\leq d~ \mbox{and}~v_j\in Q_i\}$
form a basis for $M$. Thus every $u\in M$ can be written as
\(\displaystyle \sum_{i=1}^d \langle u,\lambda_i\rangle \cdot
\mu_{i}^j\).  By (\ref{r_i}),
\(\displaystyle (\prod_{1\leq i\leq d} r_i^{\langle u,
  \lambda_i\rangle})_j=e^{u}\) for $1\leq j\leq m$ where by definition
$\lambda_i=\Lambda(Q_i)$. Thus $\bar{\phi}(\theta(e^u))=\Delta(e^u)$
so that $\bar{\phi}$ preserves the $R$-algebra structure. Hence the
lemma.\hfill$\square$

\section{A Presentation for $K_{T}^*(X)$}
Recall from \cite[Lemma $2.3$]{torusmanifold} that for each $i$,
$1 \leq i \leq d$, there exists $T$-equivariant complex line bundle
$L_i$ such that $c_1(L_i)=[V_i] \in H^2(X, \Z)$, where $[V_i]$ denotes
the cohomology class dual to $V_i$ and each $L_i$ admits an
equivariant section $\sigma_i : X \rightarrow L_i$ which vanishes
precisely along $V_i$. By \cite[Chapter VI, Theorem $2.2$]{bredon})
each $V_i$ has a closed invariant tubular neighbourhood denoted by
$N_i$ which is equivariantly diffeomorphic to the disk bundle
associated to the normal bundle $\nu_i$ of $V_i$. We denote by
$\pi : N_i \rightarrow V_i$ the projection map of the disk bundle.
Consider the trivial complex line bundle
$\epsilon:=(X \setminus \text{int}N_i) \times \C$ on
$(X\setminus \text{int}N_i)$, with the canonical $T$-action on
$(X \setminus \text{int}N_i)$ and the trivial $T$-action on the fibre
$\C$. The bundle $L_i$ is constructed by glueing $\epsilon$ and
$\pi^*(\nu_i)$ along ${\partial N_i}$. Also observe that
$N_i \setminus V_i$ does not contain any fixed points and
$L_i\mid_{V_i}\cong\nu_i$. Now recall that at a fixed point $x_j \in X^T$,
\(\displaystyle T_{x_j}(X)=\bigoplus_{ i:x_j \in V_i } \nu_i
\mid_{x_j}\) as $T$-representation space and the weight of the direct
factor $\nu_i \mid_{x_j}$ is ${\mu_i^j}$ (see the proof of \cite[Lemma
$3.5.9$]{gmukherjee}).  Since $x_j\in V_i$ if and only if $v_j\in Q_i$
we have
\begin{equation}\label{rm}[L_i \mid_{x_j}] = e^{\mu_i^j}=(r_i)_j.\end{equation}

\begin{rmk} See \cite[Section 6.1]{davisjanu} or \cite[Section
  3.2]{SU} for alternate construction of the canonical $T$-equivariant
  line bundle $L_{i}$.
\end{rmk}

\begin{thm}\label{main}
  The map
  \( \displaystyle \psi: \mathcal{K}(Q,\Lambda)\longrightarrow
  K^*_{T}(X)$ that sends $y_i\mapsto [L_i] \) for $1\leq i\leq d$ is
  an isomorphism of $RT$-algebras. Further, $\psi$ gives a
  presentation of $K_{T}^*(X)$ as a $K^*(pt)$-algebra.
\end{thm}

\noindent {\bf Proof: } Recall that the $RT$-algebra structure on
$\mathcal{K}(Q,\Lambda)$ is obtained by sending
\(\displaystyle e^u\mapsto \prod_{1\leq i\leq d} y_i^{\langle
  u,\lambda_i\rangle}\). Also the canonical $RT$-algebra structure on
$K_{T}^*(X)$ coming from the $R$-algebra structure sends $e^u\in R(T)$
to the class in $K_{T}^0(X)$ of the trivial line bundle
$L_{u}:=X\times \mathbb{C}_{u}$ on $X$ where $T$-action on fibre is
via the character $\chi_u$.  We have a canonical isomorphism of
$T$-equivariant line bundles
\(\displaystyle L_u\simeq \prod_{1\leq i\leq d}L_i^{\langle
  u,\lambda_i\rangle}\) (see \cite[(3.2)]{SU} or \cite[\mbox{arguments
  following} (3.5)]{torusmanifold}).  Since
\(\displaystyle\prod_{1\leq i\leq d} y_i^{\langle
  u,\lambda_i\rangle}\) maps to
\(\displaystyle\prod_{1\leq i\leq d}L_i^{\langle u,\lambda_i
  \rangle}\) it follows that $\psi$ is a map of $RT$-algebras.

Note that under the natural map
\(\displaystyle \iota^*: K^*_T(X) \rightarrow \prod_{i=1}^m
K_{T}^*(x_i)=R^m\) induced from the inclusion
$X^T \stackrel{\iota}{\hookrightarrow} X$, $[L_i]$ maps to
$\left( [L_i \mid_{x_j}]\right)$ which equals to $r_i$ in $\Gamma_X$
(see (\ref{rm})). Thus $\bar{\phi}=\iota^*\circ\psi$. Now, by Theorem
\ref{7}, $\bar{\phi}$ is a $R$-algebra isomorphism and hence an
$RT$-algebra isomorphism from $\mathcal{K}(Q,\Lambda)$ onto
$\Gamma_X$. Further note that under $\iota^*$, $[L_u]\in K_{T}^0(X)$
maps to
\(\displaystyle\prod_{1\leq i\leq d}r_i^{{\langle
    u,\lambda_i\rangle}}=\Delta(e^u)\). Thus $\iota^*$ is an
isomorphism of $RT$-algebras from $K_{T}^*(X)$ onto $\Gamma_X$.  It
follows that $\psi$ maps $\mathcal{K}(Q,\lambda)$ isomorphically onto
$K_{T}^*(X)$ as an $RT$-algebra.

Let $\varphi:X\longrightarrow pt$ denote the canonical
projection. Then $\varphi^*:R\longrightarrow K_{T}^*(X)$ satisfies
$\iota^*\circ \varphi^*=\Delta$ (see Remark \ref{composeconstant}).
Thus $\iota^*$ is an $R$-algebra isomorphism of $K_{T}^*(X)$ with
$\Gamma_{X}$ and hence a $K^*(pt)$-algebra isomorphism.  Since by
definition $\bar{\phi}$ is also an $K^*(pt)$-algebra isomorphism, it
follows that the isomorphism $\psi$ gives a presentation of
$K_{T}^*(X)$ as a $K^*(pt)$-algebra.  \hfill$\Box$

As a corollary we recover the $K$-ring of $X$ in terms of generators
and relations. Recall that the representation ring $RT$ acts on $\Z$
via the augmentation map $\epsilon: RT \rightarrow \Z $ taking a
virtual representation to its dimension. A $T$-space $X$ is called
weakly equivariantly formal \cite[Definition $4.1$]{HaLa}, if the map
\(\displaystyle K^*_T(X) \otimes_{RT} \Z \rightarrow K^*(X) \) induced
by the forgetful homomorphism $K^*_T(X) \rightarrow K^*(X)$ is an
isomorphism. Recall that $\mathcal{K}(Q, \Lambda)$ is a free
\(\displaystyle R=K^*(pt)\otimes_{\mathbb{Z}} RT\)-module by
\cite[Proposition $4.4$]{torusmanifold} (see also \cite[Theorem
$2.3$]{baggio2007equivariant}). Since $R$ is free as an $RT$-module it
implies in turn that $\mathcal{K}(Q, \Lambda)$ is a free
$RT$-module. Now, by Theorem \ref{main}, \(\displaystyle K^*_T(X)\) is
isomorphic to $\mathcal{K}(Q, \Lambda)$ as an $RT$-module and is hence
a free $RT$-module. Thus $Tor_{RT}^i(K^*_T(X), \Z)=0$ for all
$i>0$. It follows that the quasitoric manifold $X$ is weakly
equivariantly formal (see \cite[Proposition 4.3]{HaLa}). Denote by
$I_T:= \text{ ker}(\epsilon)$. Then the forgetful homomorphism
$K^*_T(X) \rightarrow K^*(X)$ is surjective with kernel
$I_T \cdot K_T^*(X)$ (see \cite[Proposition 4.2]{HaLa}). (See
\cite[Section 6.3]{vezzosi2003higher} for the corresponding results on
algebraic $K$-theory of smooth toric varieties).
\begin{cor}\label{maincor}(see \cite[Proposition 3.2]{SU}])
The map \(\displaystyle{\Phi}:\frac{\Z[y_1^{\pm 1}, \ldots, y_d^{\pm 1}]}{\mathcal{J}}  \longrightarrow K^*(X)\) that sends $y_i\mapsto [L_i]$ for $1\leq i\leq d$ is an isomorphism of $\mathbb{Z}$-algebras, where  where  the ideal $\mathcal{J}$ is generated by the following elements:
\begin{itemize}
\item[(i)]
  \(\displaystyle\left(1- y_{i_1} \right) \cdots \left(1- y_{i_r}
  \right), \text{ whenever } Q_{i_1} \cap \cdots \cap Q_{i_r} =
  \emptyset,\)

\item[(ii)]
  \(\displaystyle\prod_{i=1}^d y_{i}^{ \langle u, \lambda_i
    \rangle}-1 $ for $u \in M.\)
\end{itemize}
\end{cor} 

\noindent
{\bf Proof: } Note that the kernel $I_T$ of the augmentation map
$\epsilon$ is generated by $\{1 - e^u : u \in M\}$. Hence
\(\displaystyle I_T \cdot K_T^*(X)=\langle 1-\prod_{i=1}^d
y_i^{\langle u, \lambda_i \rangle} \rangle\). So \(\displaystyle K^*(X) \cong \frac{K^*_T(X)}{I_T \cdot K_T^*(X)} \cong \frac{\Z[y_1^{\pm 1}, \ldots, y_d^{\pm 1}]}{\mathcal{J}}\) and under this isomorphism $y_i$ maps to $[L_i]$.
\hfill$\Box$

\end{document}